\documentclass[12pt,reqno,twoside]{amsart}
\usepackage{amsmath}
\usepackage{amssymb}

\title{Bounded geodesics in moduli space}
\author{Dmitry Kleinbock}
\address{Brandeis University, Waltham MA
02454-9110 {\tt kleinboc@brandeis.edu}}
\author{Barak Weiss}
\address{Ben Gurion University, Be'er Sheva, Israel 84105
{\tt barakw@math.bgu.ac.il}}

\newcommand{\R}{{\mathbb{R}}}
\newcommand{\Z}{{\mathbb{Z}}}
\newcommand{\C}{{\mathbb{C}}}

\newcommand{\N}{{\mathbb{N}}}

\newcommand{\cl}{\overline}

\newcommand{\homeo}{{\operatorname{Homeo}}}

\newcommand{\Mod}{\operatorname{Mod}}
\newcommand{\SL}{\operatorname{SL}}

\newcommand{\bq}{{\mathbf{q}}}

\newcommand{\QQ}{{\mathcal Q}}
\newcommand{\QQQ}{{\til{\mathcal Q}}}

\newcommand{\LL}{{\mathcal L}}

\newcommand{\df}{{\, \stackrel{\mathrm{def}}{=}\, }}

\newcommand{\til}{\widetilde}

\newcommand{\sm}{\smallsetminus}
\newcommand{\vre}{\varepsilon}

\newcommand\hd{Hausdorff dimension}
\newcommand\nz{\smallsetminus \{0\}}
\newcommand {\equ}[1]     {\eqref{#1}}
\newcommand {\ignore}[1]  {}

\newtheorem{thm}{Theorem}

\newtheorem{cor}[thm]{Corollary}

\newtheorem{claim}{Claim}

\begin{document}


   \begin{abstract}
In the moduli space of
quadratic differentials over complex structures
on a surface,
we construct a set of full Hausdorff dimension of points
with bounded Teichm\"uller geodesic trajectories.
The main tool is quantitative nondivergence of
Teichm\"uller horocycles, due to Minsky and Weiss.
This has an application to billiards in rational polygons.

\end{abstract}

\maketitle

\section{Introduction}
Let $S$ be a compact orientable surface of genus $g$ with $n$
punctures,
where
$3g+n \geq 3$,
%
%
and let
$\QQ$ be the moduli space of unit-area holomorphic quadratic
differentials over complex structures on $S$.
This is a noncompact orbifold, on which $G \df \SL(2,\R)$ acts
continuously
with discrete stabilizers,
preserving a smooth finite measure.

For $t \in \R$, let $\displaystyle{ g_t = \left(\begin{array}{cc}
e^{t/2} & 0 \\ 0 & e^{-t/2} \end{array} \right) \in G \,, }$ and
denote the one-parameter subgroup $\{g_t : t \in \R\}$ by $F$. The
restriction of the $G$-action to $F$ defines a flow
on $\QQ$ called the {\em geodesic flow} (or sometimes the {\em
Teichm\"uller geodesic flow}). Its dynamical properties have been
extensively studied in
connection with various problems in ergodic theory and
geometry. See the surveys \cite{Howie -
handbook}, \cite{Mosher - cambridge} for more details.

It is known \cite{Masur - KC} that the
geodesic flow is ergodic
and hence a typical trajectory for the geodesic
flow is dense. It is natural (and of importance for applications) to
inquire as to the existence and abundance
of
atypical trajectories. Let us denote the
\hd\ of a subset $Y$ of a metric space $X$  by
$\dim(Y)$, and  say that $Y$ is {\it thick in} $X$
    if for any
nonempty open subset $W$ of $X$,
$
\text{dim}(W\cap
Y) = \text{dim}(W)
$
(i.e.~$Y$  has full \hd\ at any point of $X$).

In this note we discuss {\em bounded} trajectories, and prove the
following:

\begin{thm}
\label{thm: full dimension in space}
The set $\displaystyle{\{q\in \QQ : \text{the orbit }Fq\text{ is
bounded}\,\} }$ is thick in $\QQ$.
\end{thm}

Theorem \ref{thm: full dimension in space} follows from a more precise
result. Recall (see \S 2 for more details) that $\QQ$ is partitioned
into finitely many $G$-invariant suborbifolds called {\em strata}. Say
that $X \subset \QQ$ is {\em bounded in a stratum} if its closure is a
compact subset of a single stratum.

Let
$$F^+ \df \{g_t : t \geq 0\} \subset F\,.$$
For $s, \, \theta \in \R$ let
$$
h_s \df \left(\begin{array}{cc}
1 & s\\
0 & 1
\end{array}
\right), \ \ \ r_{\theta} \df
\left(\begin{array}{cc}
\cos \theta & -\sin \theta \\
\sin \theta & \cos \theta
\end{array}
\right)
\,,$$
and
$$
B\df\left\{\left(\begin{array}{cc}
a & 0\\
b & a^{-1}
\end{array}
\right):b\in\R,\ a\in\R\nz\right\}.$$

With this notation, we have:
\begin{thm}
\label{thm: more precise}
Let $q \in \QQ$. Then
\begin{itemize}
\item[(i)]
$\{s \in \R : F^+ h_sq \text{  is  bounded in  a
stratum}\, \}$ is thick in $\R$.
\item[(ii)]
For any
$f: \R \to B,$ \newline
$\{s \in \R : F^+ f(s) h_sq \text{  is  bounded in  a
stratum}\, \}$ is thick in $\R$.
\item[(iii)]
$\{\theta \in \R : F^+ r_{\theta} q \text{  is  bounded in  a
stratum}\, \}$ is thick in $\R$.
\item[(iv)]
$\{x\in G : Fxq \text{  is  bounded in  a
stratum}\, \}$ is thick in $G$.
\end{itemize}
\end{thm}

It follows from a result of Masur \cite{Masur -
log.laws} that the 
set in (iii) has measure zero, and it follows
easily that the same holds for the sets in (i), (ii), (iv).

\medskip

Note that Theorem \ref{thm: full dimension in space} follows
immediately from Theorem \ref{thm:
more precise}(iv) by considering the foliation of $\QQ$ into
$G$-orbits.
Theorem \ref{thm: more precise}(iii) has an interpretation in terms of
rational billiards,
which we now describe. Let $\mathcal{P}$ be a rational billiard table.
Fix any $\theta
\in \R$, and for any $p \in \mathcal{P}$ consider a particle moving
at a constant velocity along a ray at angle $\theta$ starting at $p$. The
law
of reflection (angle of
incidence equals angle of return)
gives rise to a dynamical system,
whose phase space $\mathcal{P}_{\theta}$ is
   the complement of a countable
set of lines in finitely many copies of $\mathcal{P}$.
See \cite{Howie - handbook} for more details.
Let $d$
denote the restriction of the Euclidean metric to
$\mathcal{P}_{\theta}$, and let
$b^\theta_t:
   \mathcal{P}_{\theta} \to
\mathcal{P}_{\theta}$ denote the time $t$ map of this flow.


It can be shown (see e.g.\ \cite[\S 2]{Boshernitzan - rank two}) that
every billiard trajectory is {\em recurrent\/}, that is,
  $\liminf_{t\to\infty}d(p,
b^\theta_tp) = 0$ for any  $p \in
\mathcal{P}_{\theta}$. Furthermore, one can use  \cite[Theorem
1.5]{Boshernitzan - recurrence} and the existence of one-dimensional
Poincar\'e sections for the flow
$(\mathcal{P}_{\theta},b^{\theta}_t)$ (see a discussion in \cite[\S
1.7]{Howie - handbook}) to derive a quantitative version  for generic
trajectories: namely,
given $\mathcal{P}$ and $\theta$, there is a constant $C$, such that
for almost every $p \in
\mathcal{P}_{\theta}$ one has
$$\liminf_{t\to\infty}t\cdot d(p, b^\theta_tp)
\le C
\,.$$

On the other hand it turns out that for a thick set of angles
$\theta$ it is not possible to replace $C$ in the right hand
side of the above inequality by an arbitrary small number:

\begin{cor}
\label{cor: billiards}


The set
$$\{\theta : \exists\, c > 0 \text{ such that
}\liminf_{t\to\infty}\,{{t}}\cdot d(p, b^\theta_tp)  \ge c\quad
\forall\,p \in
\mathcal{P}_{\theta}\}$$ is thick in $\R$.

\end{cor}


Note that if $\mathcal{P}$ is rectangular (or more generally
{\em integrable}, see \cite{Howie - handbook}), then $b^{\theta}_t:
\mathcal{P}_{\theta} \to \mathcal{P}_{\theta}$ is isomorphic to the
flow along a ray in the torus $\R^2/\Z^2$, and Corollary \ref{cor:
billiards} 
is equivalent to the well-known fact (see \cite{Schmidt} and
the  references therein) that
the set of
{\em badly approximable\/} $\theta
\in \R$  (that is, those for which there exists
$c>0$ with $\displaystyle{|\theta - \frac{p}{q}| \geq
\frac{c}{q^2}}$ for all integers $p, \,q$)
is thick.

\medskip

Our work is  part of a fruitful interaction between the study of
dynamics on quadratic differential spaces and dynamics on homogeneous
spaces of Lie groups, see the survey \cite{Cambridge} and the references
therein. Bounded trajectories on homogeneous spaces were extensively
studied due to their relation to linear forms which are badly
approximable by rationals, and results analogous to ours were obtained
by W.\ Schmidt \cite{Schmidt}, S.G.\ Dani \cite{Dani}, and, in most
general form, in
joint work of the first-named author and G.A.\ Margulis \cite{Dima -
Margulis}. See  \cite[\S 4.1b]{KSS survey} for references to other
related results.

Our construction of bounded trajectories is similar in some respects
to the one used in
\cite{Dima - Margulis}, but whereas in that paper a crucial estimate
was obtained by using the mixing property of the flow, in this paper
the estimate follows from quantitative nondivergence results for
horocycles obtained in \cite{with Yair}. This
approach,
powered by nondivergence estimates for unipotent
trajectories
of homogeneous flows \cite{KM - manifolds}, can also be used to
reprove the abundance of bounded orbits for
partially
hyperbolic diagonalizable
flows on
$\SL(n,\R)/\SL(n,\R)$,
and to obtain new applications to number theory.
It was originally
introduced in the talk given by the first-named author  in July 2000 at a
Euroconference on  Ergodic Theory, Number Theory and
Geometric Rigidity  (Issac Newton Institute, Cambridge, UK).
%

\medskip

{\bf Acknowledgements:}
We thank Lee Mosher for getting us
interested in this problem, and for encouraging remarks. We also thank
Howard Masur and Yair Minsky for teaching us about quadratic
differentials. This research was supported by BSF grant 2000247.

\medskip

\section{Terminology}
Let $\Mod(S)$ be the mapping class group, let $\til{\QQ}$ be the space
of unit-area holomorphic quadratic differentials over complex
structures on $S$, so that $\QQ = \QQQ/\Mod(S)$.
Let $\pi: \til{\QQ}
\to \QQ$ be the natural quotient map. The space $\QQQ$ is a manifold
on which $G$ acts continuously, and $\pi$ is $G$-equivariant. Each
$\bq \in \QQQ$ can be defined as an equivalence
class of atlases of charts of the following type. Outside a
finite singularity set $\Sigma = \Sigma(\bq) \subset S$ we have charts
$S \sm \Sigma \to \C$ for which the transition functions $\C \to \C$
are of the form $z \mapsto \pm z +c$. The charts are required to
have $k$-pronged singularities around points of $\Sigma$, with $k \geq
3.$ At the punctures they are required to extend to the `filled-in'
surface, and the punctures may have $k$-prongs with $k \geq 1$.
The equivalence
relation is given by the natural action of
$\homeo_0(S)$ by pre-composition on each chart. An atlas is said to be
orientable if the transition functions can be taken to be of the form
$z \mapsto z +c$.
    See \cite[\S 4]{Howie -
handbook} or \cite[\S 4]{with Yair} for details.

The data consisting of the number of singularities of each type, and
the orientability of the atlas, assumes one of finitely many
possibilities. A level set for this data is called a {\em
stratum} in $\QQQ$. If $\mathcal{M} \subset \QQQ$ is a stratum we also
call its projection $\pi(\mathcal{M}) \subset \QQ$ a {\em stratum} in
$\QQ$. Any stratum in $\QQQ$ is a submanifold and any
stratum in $\QQ$ is a sub-orbifold. The strata
are invariant under the $\SL(2,\R)$-action.

A {\em   saddle connection\/}
for
$\bq$ is a path
$\delta:(0,1)\to S \sm \Sigma$ whose image in any chart is a Euclidean
straight line and which extends continuously to a map from $[0,1]$ to the
      completion of $S$, mapping $\{0,1\}$ to singularities or punctures.
In all cases except $g=1, \, n=0$, there are saddle connections on the
surface. Since our results are immediate for this case, we assume from
now on that $(g,n) \neq (1,0)$.

     Integrating the local projections of
$d\delta$
one obtains
a vector, denoted by
$\big(x(\delta,\bq),y(\delta,\bq)\big)$, well-defined up to a multiple of
$\pm 1$. Let $\LL_{\bq}$ denote the set of all saddle connections for
$\bq$. Since linear maps do not change the property of being a
straight line, we may identify $\LL_\bq$ with $\LL_{A\bq}$ for any
$A\in G$.  The action of
$G$ transforms the components of saddle connections linearly:
\begin{equation}\label{eq: action on coordinates}
\text{for }A \in G\,, \ \left(
\begin{array}{c}
x(\delta,A\bq)\\
y(\delta,A\bq)
\end{array}
\right)
=\pm A
\left(
\begin{array}{c}
x(\delta,\bq)\\
y(\delta,\bq)
\end{array}
\right)\,.
\end{equation}

For a saddle connection $\delta\in\LL_\bq$, define its `length'
$l(\delta,\bq)$ by
\begin{equation*}
\label{eq: length def}
l(\delta,\bq)  \df \max \big(|x(\delta,\bq)|,|y(\delta,\bq)|\big)\,.
\end{equation*}

For $\vre>0$, let
$$K_{\vre} \df \pi \left(\{\bq \in \til{\QQ}:
\, \forall \delta \in \LL_{\bq}, \, l(\delta, \bq) \geq \vre\} \right).
$$
It is known that each $K_{\vre}$ is compact.
Moreover, for each
stratum $\mathcal{M}$ in $\QQ$, each compact $K \subset \mathcal{M}$
is contained in $K_{\vre}$ for some $\vre$, and
$$\cl{\pi(\mathcal{M}) \cap K_{\vre}} \subset \pi(\mathcal{M}).$$
    Thus $X \subset \mathcal{M}$ is bounded in a stratum if and only if
$X \subset K_{\vre}$ for some $\vre>0$.

\section{From nondivergent horocycles to bounded geodesics}

Recall that the restriction of the $G$-action to the subgroup $\{h_s
:s \in \R\}$ defines the so-called {\em horocycle flow}. It was proved
in \cite{Veech}, using arguments of \cite{KMS}, that this flow does not
admit divergent trajectories. The main
ingredient for the proof of Theorem \ref{thm: more precise}
is a quantitative strengthening of this result, along the lines of
\cite{KM - manifolds}.
   Here $|\cdot |$ stands for Lebesgue measure on $\R$.

\begin{thm}{\cite[Thm. 6.3]{with Yair}}
\label{thm: nondivergence}
There are positive constants $C, \, \alpha, \,
\rho_0$, depending only on $S$, such that if $\bq \in \QQQ$, an
interval $I \subset \R$, and $0<\rho \leq \rho_0$ satisfy:
\begin{equation}
\label{eq-condition}
{\rm for \ any \ } \delta \in \LL_{\bq}, \ \sup_{s \in I} l(h_s\bq, \delta)
\geq \rho,
\end{equation}
then for any $0<\vre<\rho$ we have:
\begin{equation}
\label{eq - precise yet long}
|\{s \in I: h_s \pi(\bq) \notin K_{\vre} \}| \leq
C\left(\frac{\vre}{\rho}\right)^{\alpha}|I|.
\end{equation}

\end{thm}

\medskip

{\bf Proof of Theorem \ref{thm: more precise}.}
We start with part (i).
Let $J \subset \R$ be an interval and let $q_0 = \pi(\bq_0) \in \QQ$. We
need to prove that the Hausdorff dimension of the set
\begin{equation}
      \label{eq: set}
\{ s \in J : \, F_+ h_{s}q_0 \mathrm{\ is \
bounded \ in \ a \ stratum} \,\}
\end{equation}
is equal to $1$.
Denote the length of $J$ by $r$, and assume, as we may without loss of
generality, that $J = [0,r]$ and $r\le 1$.

Let $\eta>0$ be arbitrary, and
let $C, \, \alpha, \, \rho_0$ be as in Theorem \ref{thm:
nondivergence}. Choose $\vre>0$ so that
$q_0\in K_{\vre}$, and also
$$C \left (\frac{\vre}{ \rho_0}\right) ^{\alpha} < \eta\,,$$
and let
\[
K' \df K_{\vre}, \ \ \ \ \ K''
\df
\bigcup_{0\le s \le r
} h_s K'\,.
\]

Then there is $\vre_0$, depending only on $\vre$,
such that $K'' \subset K_{\vre_0}$. Hence for
each $q \in K''$, each $\bq \in \pi^{-1}(q),$ and
each $\delta  \in
\LL_{\bq}$, we have
\begin{equation}
\label{eq: vre_0}
l(\delta , \bq) = \max \big(|x(\delta,\bq)|,|y(\delta,\bq)|\big) \geq
\vre_0\,.
\end{equation}

\begin{claim}
\label{claim: achieve large}
There exists $t_1>0,$ depending only on
$\vre_0$ and $r$,
     such that for
every
$q \in K'', \, \bq \in \pi^{-1}(q)$ and $\delta \in
\LL_{\bq}$,
\begin{equation}
\label{eq: ii}
     \sup_{0\le s \le r}
|x(\delta,
g_{t_1} h_{s}  \bq)|
\geq \rho_0\,.
\end{equation}

      \end{claim}

%

To prove the claim,
let us first show that
\begin{equation}
\label{eq: i}
\max\big(|x(\delta, \bq)|,|x(\delta, h_{r} \bq)|\big) \geq
\vre_0r/2\,.
\end{equation}
Let
$x=x(\delta,  \bq)$, $y=y(\delta,
     \bq)$. If $|x| \geq
\vre_0r/2$, there is nothing to prove. So by (\ref{eq: vre_0}) we may
assume that $|y| \geq \vre_0$. But then
\[
     |x(\delta, h_{r} \bq)|
\stackrel{\equ{eq: action on
coordinates}}{=} |x   + ry| \ge \vre_0r -
\vre_0r/2=\vre_0r/2 \,.
\]
Setting $t_1 \df 2\log \frac{2\rho_0}{\vre_0r}$, one obtains
(\ref{eq: ii}) from (\ref{eq: i}) via (\ref{eq: action on
coordinates}). This proves the claim.

\medskip

Enlarging $t_1$ if necessary, assume that
$$N \df e^{t_1} \in \N\,.$$

\begin{claim}
\label{claim: return}
Suppose $q' \in \QQ
$ and $t_0 \geq 0$ are
such that
$g_{t_0} h_s q' \subset K''$ for any $\,0 \le s \le re^{-t_0}$.
Subdivide
$[0,re^{-t_0}]$
into $N$
subintervals $J_1, \ldots, J_N$ of length $re^{-t_0}/N=re^{-(t_0+t_1)}$,
and let
$$\mathcal{I} \df \left \{i \in \{1, \ldots, N\} :
g_{t_1+t_0}h_{s}q' \in K'' \ \forall\,s\in J_i\right \}
.$$
Then $\# \mathcal{I} \geq (1 -\eta)N.$
\end{claim}

In order to prove this claim,
we apply Theorem \ref{thm: nondivergence} with
$I = [0,re^{t_1}]$ and
$q =  g_{t_1+t_0}q'$, that is, to the horocycle
\[
\{h_s g_{t_1+t_0}q' : 0 \le s \le re^{t_1}\} = \{g_{t_1}h_s g_{t_0}q':
0 \le s
\le r\}\,.
\]
     Assumption (\ref{eq-condition}), with
$\rho=\rho_0$, is valid by (\ref{eq: ii}).

We obtain:
\begin{equation}
      \label{eq: intersect most of the time}
\begin{split}
|\{0 \le s \le re^{-t_0}:   g_{t_1+t_0}h_sq' \notin K' \}| &=
e^{-(t_1 + t_0)}|\{0 \le s \le re^{t_1}:  h_s g_{t_1+t_0}q' \notin K'
\}|\\
     &\leq
e^{-(t_1 + t_0)} C\left(\frac{\vre}{\rho_0}\right)^{\alpha}re^{t_1} < \eta
re^{-t_0}\,.
\end{split}
\end{equation}

Now let
$$\mathcal{I}_0 \df \{ i \in \{1, \ldots, N\} : \,
g_{t_1+t_0} h_{s}q' \in K' \text{ for some }s\in J_i \}\,.
$$
Then (\ref{eq: intersect most of the time}) implies that
$\# \mathcal{I}_0 \geq (1-\eta)N.$
A simple calculation
using
the definition of $K''$ shows that
$\mathcal{I}_0 \subset \mathcal{I}$, and the claim follows.

\medskip

We now construct a subset $D \subset J$, of large Hausdorff dimension,
by iterating Claim \ref{claim:
return}. At stage $m$ of the construction we have a subdivision of $J$
into $N^m$ subintervals of length $N^{-m}$, and a collection of
$[(1-\eta)N]^m$ of these subintervals with the property that if $I$ is
one of these intervals then
$$ g_{kt_1}h_s q_0 \subset  K'' \text{ for any }s\in I\text{ and }\ k=0,
\ldots, m\,.$$ For each subinterval $I$ belonging to this collection, we
apply Claim
\ref{claim: return} to $I$ shifted by its left endpoint, with $t_0=mt_1$
and an appropriate
$q' \in \pi(\mathcal{M})$, and obtain that there are at least $(1-\eta)N$
subintervals of $I$ which satisfy the conclusion of the claim. The union
of these subintervals over all subintervals $I$ comprises the collection
for the \linebreak
$(m+1)$-st step.

Let $D$
be the set of intersection points of all
sequences of nested subintervals in the above construction.

A well-known result
     (see  \cite[Proposition 2.2]{McMullen})
implies that the Hausdorff dimension of $D$ is at
least
$$1-\frac{\log(1-\eta)}{-\log(N)} \geq 1+\log(1-\eta)\,.$$
If $s \in
D$, then
$g_{kt_1}h_{s}q \in K''$ for any $k \in \Z_+$.
This
implies that
$$F_+ h_{s} q_0 \subset K \df \bigcup_{0\le t \le t_1
} g_t K''\,,$$
where $K$ is compact.

Since $\eta$ was chosen arbitrarily and since $1+ \log(1-\eta)$ tends
to $1$ as $\eta\to 0$, we obtain that the set \equ{eq: set} has Hausdorff
dimension $1$, as required.

\medskip

For part (ii), note that for every $b \in B$, the set
$$\{g_tbg_{-t} : t \geq 0\} \subset G$$
is precompact, hence
$$F^+ f(s) h_sq = \big\{\big(g_t f(s)
g_{-t}\big) g_t  h_sq : t\ge 0\big\}$$ is bounded if and only if so is
$F^+ h_sq$; thus (ii) follows from (i). The claim (iii) is a special
case of (ii): it is easy to represent $r_\theta$ in the form $r_\theta =
f(\theta) h_\theta$ where $f(\theta) \in B$.

Finally, to prove (iv), we
follow the argument of \cite[\S 1.5]{Dima -
Margulis}. Take
$U\subset G$  of the form
$V^- V V^0$, where $V^-$, $V$  and $V^0$ are
neighborhoods of identity in
$$
H^-\df \left\{\left(\begin{array}{cc}
1 & 0\\
s & 1
\end{array}
\right) : s\in \R \right\}
, \ \
H \df \{h_s: s \in \R\} \ \ \mathrm{and \ } F \
\mathrm{respectively},$$ such
that the multiplication map $V^-\times V\times V^0\to U$
is bi-Lipschitz. It is enough to show that for any $q$,
\begin{equation}
      \label{eq: dim 3}
\dim\big(\{h \in U : Fhq\text{ is
bounded}\,\}\big) = 3\,.
\end{equation}
Let $C$ be the set $\{h\in U : F^+hq\text{ is
bounded}\,\}$. From part (i) it
follows that for any $h^0\in V^0$, $\dim(C\cap Vh^0) = 1$. In view of the
slicing properties of the \hd\ (see \cite{Falconer} or \cite[Lemma
1.4]{Dima - Margulis}),  one deduces that
\begin{equation}
      \label{eq: dim 2}
\dim(C\cap V V^0) = 2\,.
\end{equation}

Now for any $h\in C$ choose a neighborhood $V^-(h)$
of identity in $H^-$ such that
$V^-(h)h\subset U$. Then, by the argument used in the proof of (ii), one
necessarily has
$V^-(h)h\subset C$. Applying part (i) of the theorem
with
$$F^-\df \{g_t : t\le 0\}$$ in place of $F^+$ and $H^-$ in place
of $H$ (this can be justified e.g.\ by considering an automorphism
$q\mapsto \left(\begin{array}{cc}
0 & 1\\
1 & 0
\end{array}
\right)q$ of $\QQ$, which leaves $G$-orbits and lengths of saddle
connections invariant, and sends the action of $H$ and $F^+$ to that of
$H^-$ and $F^-$ respectively),
one concludes that
$$
\forall\,h\in C\cap V V^0\quad\dim\big(\{h^- \in V^- : Fh^-
hx\text{ is bounded}\,\}\big) = 1\,,
$$
which, together with \equ{eq: dim 2}, yields \equ{eq: dim 3}.
\qed

%
%

\section{Application to rational billiards}


Recall that any $\bq \in \QQQ$ defines a pair of transverse measured
foliations on $S \sm \Sigma(\bq)$, called the {\em horizontal\/} and
{\em vertical\/} foliations respectively, and that $g_t$ acts by
contracting (resp.\ expanding) the measure transverse to horizontal
(resp.\ vertical) leaves by a factor of $e^{t/2}$. The Euclidean metric
on
   a billiard table
$\mathcal{P}$ determines a flat Riemannian metric on
$S
\sm
\Sigma (\bq)$. The length of $\delta \in \LL_{\bq}$ with respect to this
metric is at least as large as $l(\delta, \bq)$. Recall also that any
line segment is a length-minimizing path in its homotopy class with
endpoints fixed, and any free homotopy class has a length-minimizing
representative which is a concatenation of saddle connections.

Given $\mathcal{P}$, there is a standard construction (see
\cite{KMS} or \cite{Howie - handbook}) of a
closed surface $S$ of genus $g \geq
1$, with a  quadratic differential $\bq$, such
that, for each $\theta \in
\R$, the  flow $b^{\theta}_t : \mathcal{P}_{\theta} \to
\mathcal{P}_{\theta}$ is isomorphic to the unit-speed flow in $S \sm
\Sigma(\bq)$ along vertical leaves for $r_{\theta} \bq$.
The case $g=1$ occurs if and only if $\mathcal{P}$ is integrable, and
in this case Corollary \ref{cor: billiards} follows
from the fact that badly approximable numbers form a thick subset of
the reals. So assume $g \geq 2$.



\medskip

{\bf Proof of Corollary \ref{cor: billiards}.} Suppose $\theta$
belongs to the thick subset of $\R$ for which $F^+ r_{\theta}
\pi(\bq)$ is bounded in a stratum. Then there exists
$0<\vre<1$ with
\begin{equation}
\label{eq: fix the compact}
F^+r_{\theta}\pi(\bq) \subset K_{\vre}\,.
\end{equation}
Let $c \df {{\vre^2}/{2}<1}$; we claim that
${d(b^\theta_tp,p)\ge{c}/{t}}$ for any $p
\in
\mathcal{P}_{\theta}$  and $t\geq 1$.
Indeed, assuming  that there are $p \in \mathcal{P}_{\theta}$ and $t\geq
1$ with
${d(b^\theta_tp,p) < {c}/{t}}$, let
$$t_0 \df 2 \log t - \log c \geq 0\,.$$

With respect to $r_{\theta} \bq$, let $\mathcal{C}_1$ be the vertical
segment from $p$ to $b^\theta_tp$, let $\mathcal{C}_2$ be a
length-minimizing path from $b^\theta_tp$ to $p$, and let $\mathcal{C}$
be the concatenation of
$\mathcal{C}_1$ and $\mathcal{C}_2$. The length of $\mathcal{C}_1$
with respect to $r_{\theta} \bq$ is $t \geq 1$ and the length of
$\mathcal{C}_2$ is less
than $c/t<1$. Since $\mathcal{C}_1$ is length-minimizing in its
homotopy class with fixed endpoints, $\mathcal{C}$ is homotopically
nontrivial.

With respect to $g_{t_0} r_{\theta} \bq$, the length of
$\mathcal{C}_1$ is $te^{-t_0/2} = \sqrt{c}$, and the length of
$\mathcal{C}_2$ is less than
$e^{t_0/2}c/t=\sqrt{c}$. Thus the length of $\mathcal{C}$ with respect
to $g_{t_0}r_{\theta} \bq$ is less than $\vre$. Since
$\mathcal{C}$ is not homotopically trivial, the shortest
representative in its free homotopy class is a nontrivial
concatenation of saddle
connections. In particular there is $\delta \in
\LL_{\bq}$ with $l(\delta, g_{t_0}r_{\theta} \bq) < \vre$, which
contradicts (\ref{eq: fix the compact}). \qed




\end{document}